\UseRawInputEncoding
\documentclass[11pt]{amsart}
\usepackage{mathrsfs}
\usepackage{amssymb}

\usepackage{amscd}
\usepackage{amsmath, amssymb}
\usepackage{amsfonts}
\usepackage[colorlinks,linkcolor=blue,citecolor=blue, pdfstartview=FitH]{hyperref}
\usepackage[all]{xy}

\numberwithin{equation}{section}

\theoremstyle{plain}
\newtheorem{thm}{Theorem}[section]
\newtheorem{theorem}[thm]{Theorem}
\newtheorem{lemma}[thm]{Lemma}

\newtheorem{proposition}[thm]{Proposition}

\theoremstyle{definition}

\newtheorem{remark}[thm]{Remark}

\newtheorem{definition}[thm]{Definition}

\newtheorem{defn-thm}[thm]{Definition-Theorem}

\begin{document}

\title{On almost nonpositive $k$-Ricci curvature}
\author{Kai Tang}
\address{Kai Tang. College of Mathematics and Computer Science, Zhejiang Normal University, Jinhua, Zhejiang, 321004, China} \email{{kaitang001@zjnu.edu.cn}}
\thanks{\textbf{Foundation item:} Supported by National Natural Science
Foundation of China (Grant No.12001490).}

\maketitle

\markleft{On almost nonpositive $k$-Ricci curvature} \markright{On almost nonpositive $k$-Ricci curvature}

\begin{abstract}
Motivated by the recent work of Chu-Lee-Tam on the nefness of canonical line bundle for compact K\"{a}hler manifolds with nonpositive $k$-Ricci curvature, we consider a natural notion of {\em almost nonpositive $k$-Ricci curvature},
which is weaker than the existence of a K\"{a}hler metric with nonpositive $k$-Ricci curvature. When $k=1$, this is just the {\em almost nonpositive holomorphic sectional curvature} introduced by Zhang. We firstly give a lower bound for  the existence time of the twisted K\"{a}hler-Ricci flow when there exists a K\"{a}hler metric with $k$-Ricci curvature bounded from above by a positive constant. As an application, we prove that a compact K\"{a}hler manifold of almost nonpositive $k$-Ricci curvature must have nef canonical line bundle.
\end{abstract}
\maketitle

\section{Introduction and statement of result}
In an attempt to generalize the hyperbolicity of Kobayashi to the $k$-hyperbolicity, Ni \cite{N1} introduced the concept of {\em $k$-Ricci curvature}. Given a compact K\"{a}hler manifold $(M^{n}, h)$ with K\"{a}hler form $\omega$ and Chern curvature tensor $R$. The $k$-Ricci curvature $Ric_{k}$ ($1\leq k\leq n$) is defined as the Ricci curvature of the $k$-dimensional holomorphic subspaces of the holomorphic tangent bundle $T^{\prime}M$. Clearly, $Ric_{1}$ is just the holomorphic sectional curvature $H(X)$ and $Ric_{n}$ conincides with the Ricci curvature Ric of $M$. Hitchin \cite{Hit} showed that $H$ and Ric are independent to each other by an example.

There are many imortant results on a compact K\"{a}hler manifold with $Ric_{k}\geq 0$ or $Ric_{k}\leq 0$. It was proved by Yang \cite{Yang1} that a compact K\"{a}hler manifold with $Ric_{1}>0$ must be projective and rationally connected, confirming a conjecture of Yau \cite{Yau1}. In \cite{N2}, Ni showed that it is also true if $Ric_{k}>0$ for some $1\leq k\leq n$. In their recent breakthrough \cite{WY1}, Wu-Yau confirmed a conjecture of Yau that a projective K\"{a}hler manifold with $Ric_{1}<0$ must have ample canonical line bundle. Tosatti-Yang \cite{TY} was able to drop the projectivity assumption in Wu-Yau theorem. They also proved that a compact K\"{a}hler manifold with $Ric_{1}\leq 0$ must have nef canonical bundle. In a recent preprint\cite{CLT}, Chu-Lee-Tam proved that a compact K\"{a}hler manifold with $Ric_{k}<0$ ($Ric\leq0$) have ample (nef) canonical bundle. Li-Ni-Zhu \cite{LNZ} also gave an alternate proof to the results of Chu-Lee-Tam. For more related works, we refer readers to \cite{HLW,HW,DT,WY2,Yang2,YZ,Nom,Tang,LS,NZ1,NZ2}.

As we all know, the ampleness of the canonical line bundle is equivalent to the existence of one K\"{a}hler metric with negative Ricci curvature and  the nefness of the canonical bundle is defined by taking \textbf{limits of a family of K\"{a}hler classes.} Therefore, it may not be ``best" choice to imply the nefness of the canonical line bundle from the nonpositivity of the $k$-Ricci curvature of \textbf{one K\"{a}hler metric.} A possible natural question is whether there is a condition in terms of $k$-Ricci curvature which is weaker than the existence of a K\"{a}hler metric with $Ric_{k}\leq 0$, but it can guarantee the nefness of the canonical line bundle. For the case of $Ric_{1}\leq 0$, Zhang \cite{Zhang1} was the first to consider this problem and he defined the concept of {\em almost nonpositive holomorphic sectional curvature} (namely, almost nonpositive $1$-Ricci curvature). He also proved that a compact K\"{a}hler manifold of almost nonpositive holomorphic sectional curvature has a nef canonical line bundle. Zhang-Zheng \cite{ZZ} also studied the compact K\"{a}hler manifolds with almost quasi-negative holomorphic sectional curvature.
Motivated by the  work of Zhang \cite{Zhang1}, we introduce a natural notion of {\em almost nonpositive $k$-Ricci curvature}.

Let us  recall the concept of $k$-Ricci curvature on a compact K\"{a}hler manifold $(M^{n}, h)$ introduced in \cite{N1}. We denote Chern curvature tensor as $R$. For a point $p\in M$, let $\Sigma\in T^{\prime}_{p}M$ be a $k$-dimensional subspace.
The $k$-Ricci curvature of the K\"{a}hler metric $h$ on $\Sigma$ is
\begin{align}
Ric_{k}^{h}(p,\Sigma)(X,\overline{Y})=tr_{h}R(X,\overline{Y},\cdot,\cdot) \,\,\,,\nonumber
 \end{align}
for $X,Y\in\Sigma$ where the trace is taken with respect to $h|_{\Sigma}$. For any $k$-dimensional subspace $\Sigma\in T^{\prime}_{p}M$ at any point $p$, $X\in\Sigma \backslash\{0\}$. If the following inequality
\begin{align}
\frac{Ric_{k}^{h}(p,\Sigma)(X,\overline{X})}{|X|^{2}_{h}}\leq \lambda \,\,\,\, (\geq \lambda)   \nonumber
 \end{align}
holds, we denote it as $Ric^{h}_{k}\leq \lambda$ ($\geq \lambda$). we set
\begin{align}
\mu_{h}(k)=\sup_{x\in M}\{\sup Ric^{h}_{k}|_{x}\} \,\,\,,\nonumber
 \end{align}
where it means the maximal value of the $k$-Ricci curvature of $h$ on the compact manifold $M$.

\begin{definition}\label{1.1} Let $(M^{n}, \omega_{0})$ be a compact K\"{a}hler manifold.
\begin{itemize}
\item[(a)] Let $[\omega]$ be a K\"{a}hler class on $M$. We define the number $\mu_{[\omega]}(k)$ as follows:
\begin{equation}
\mu_{[\omega]}(k)=\inf\{\mu_{\omega^{\prime}}(k)| \text{$\omega^{\prime}$ is a K\"{a}hler metric in $[\omega]$} \} \,\,\,.\nonumber
 \end{equation}
\item[(b)] If there exist a sequence number $\varepsilon_{i}\searrow 0$ and a sequence of K\"{a}hler class $\alpha_{i}$ on $M$ such that $\mu_{\alpha_{i}}(k)\alpha_{i}<\varepsilon_{i}[\omega_{0}]$, we say
that $M$ \textbf{is of almost nonpositive $k$-Ricci curvature}.

\item[(c)] If $\mu_{[\omega]}(k)=0$, we say that \textbf{the K\"{a}hler class $[\omega]$ is of almost nonpositive $k$-Ricci curvature}.
\end{itemize}
\end{definition}

In this note, we first prove the following property, which shows that $\mu_{[\omega]}(k)$ is well defined.

\begin{proposition}\label{1.2} Let $(M^{n}, \omega_{0})$ be a compact K\"{a}hler manifold. Let $\alpha$ be any K\"{a}hler class  on $M$. For $1\leq k\leq n$, we have $\mu_{\alpha}(k)>-\infty$.
\end{proposition}

\begin{remark}
(1) Obviously, if $\mu_{\alpha}(k)=0$, the $M$ must be of almost nonpositive $k$-Ricci curvature which means that the definition (b) is weaker than the definition (c). $\mu_{\alpha}(k)<0$ if and only if there exists a K\"{a}hler metric $\omega^{\prime}\in \alpha$ of negative $k$-Ricci curvature.
Clearly, the condition $\mu_{\alpha}(k)=0$ is by
definition a condition weaker than the existence of a K\"{a}hler metric $\omega^{\prime}\in \alpha$ of nonpositive $k$-Ricci curvature. \\(2) We know that nefness
is a positivity at the level of $(1,1)$-classes, not $(1, 1)$-forms, so the definition of $(b)$ seems reasonable where it is also a definition at  the level of $(1,1)$-classes. When $k=1$, it is just the concept of {\em almost nonpositive holomorphic sectional curvature}. In \cite{Zhang1}, Zhang gave
a lot of good properties and applications of the concept.
\end{remark}

In \cite{CLT}, Chu-Lee-Tam used the twisted K\"{a}hler-Ricci flow to study the compact K\"{a}hler manifolds with nonpositive $k$-Ricci curvature. By employing this method, we obtain the following result.
\begin{theorem}\label{1.3} A compact K\"{a}hler manifold $M^{n}$ of almost nonpositive $k$-Ricci curvature must have nef canonical line bundle.
\end{theorem}

In particular, from the above Theorem \ref{1.3}, we easily get the following result:
\begin{theorem} A compact K\"{a}hler manifold $M^{n}$ admitting a K\"{a}hler class $\alpha$ of almost nonpositive $k$-Ricci curvature must have nef canonical line bundle.
\end{theorem}

\begin{remark}\label{1.31}
Note that when $k=1$, the above results was proved by Zhang \cite{Zhang1}. Theorem \ref{1.3} is also generalization of  Chu-Lee-Tam's result \cite{CLT}. For $k=n$, $Ric_{n}$ is just the Chern Ricci curvature. If a K\"{a}hler class $\alpha$ is of almost nonpositive Ricci curvature, then there exists a sequence of K\"{a}hler metric $\omega_{\varepsilon}\in \alpha$ such that $Ric^{\omega_{\varepsilon}}<\varepsilon\omega_{\varepsilon}$ for any $\varepsilon>0$, namely, $2\pi c_{1}(K_{M})+\varepsilon\alpha>0$. So the canonical line bundle must be nef.
\end{remark}

To see Theorem \ref{1.3}, we give a useful proposition on a lower bound of the twisted K\"{a}hler-Ricci flow, which might have other applications.
\begin{proposition}\label{1.4} Let $(M^{n}, \omega_{0})$ be a compact K\"{a}hler manifold and $\widehat{\omega}$ be a K\"{a}hler metric of the K\"{a}hler class $[\omega_{0}]$. For a fixed integer $k$ with $1<k< n$,
we assume that $A$ is  the maximal value of the $k$-Ricci curvature of $\widehat{\omega}$ on $M$. Set $A>0$. If there exists a positive constant $\delta>0$, such that $\delta[\omega_{0}]+2\pi c_{1}(K_{M})>0$. Then we
can find a function $\upsilon\in C^{\infty}(M)$ which satisfies $\delta\widehat{\omega}-Ric(\widehat{\omega})+\sqrt{-1}\partial\overline{\partial}\upsilon>0$
 and $\eta=\frac{k-1}{2(n-k)}\sqrt{-1}\partial\overline{\partial}\upsilon$, such that the twisted K\"{a}hler-Ricci flow running from $\widehat{\omega}$,
\begin{equation}\label{1.5}
\left\{\begin{split}
&\partial_{t}\omega(t)=-Ric(\omega(t))-\eta \\
&\omega(0)=\widehat{\omega}\,\,\,,
\end{split}\right.
 \end{equation}
 exists a smooth solution on $M\times [0, \frac{2n(k-1)}{(n-1)[2nA+(k-1)\delta]})$.
\end{proposition}

\section{Proof of Proposition \ref{1.2}}
Before proving Proposition \ref{1.2}, we give some algebraic estimates which are proved by Chu-Lee-Tam \cite{CLT}.
They are useful in obtaining key estimates for the twisted K\"{a}hler-Ricci flow.
\begin{lemma}\cite{CLT}
Let $(M^{n}, h)$ be a compact K\"{a}hler manifold with $Ric_{k}(X,\overline{X})\leq-(k+1)\sigma|X|^{2}$, $\sigma\in \mathbb{R}$. Then the following inequality holds
\begin{equation}\label{2.1}
(k-1)|X|_{h}^{2}Ric(X,\overline{X})+(n-k)R(X,\overline{X},X,\overline{X})\leq-(n-1)(k+1)\sigma|X|_{h}^{4}.
 \end{equation}
\end{lemma}
Furthermore, by using the Royden's trick \cite{Roy}, the following result holds.
\begin{lemma}\cite{CLT} \label{2}
Let $(M^{n}, h)$ be a compact K\"{a}hler manifold with $Ric_{k}(X,\overline{X})\leq-(k+1)\sigma|X|^{2}$, $\sigma\in \mathbb{R}$. If $g$ is another K\"{a}hler metric, then the following inequality holds
\begin{align}\label{2.2}
2g^{i\overline{j}}g^{k\overline{l}}R_{i\overline{j}k\overline{l}}\leq &\frac{-(n-1)(k+1)
\sigma}{n-k}((tr_{g}h)^{2}+|h|_{g}^{2}) \\
&-\frac{k-1}{n-k}(tr_{g}h)\cdot(tr_{g}Ric)-\frac{k-1}{n-k}\langle h, Ric\rangle_{g} \,\,\,, \nonumber
 \end{align}
 where R is the Chern curvature tensor of $h$ and Ric is the Ricci curvature of $h$.
\end{lemma}
If the above Lemma \ref{2} satisfies $g=h$, it can also imply the relation on $Ric$ and scalar curvature $S$ under the assumption $Ric_{k}(X,\overline{X})\leq-(k+1)\sigma|X|^{2}$.
\begin{lemma}\cite{CLT} Let $(M^{n}, h)$ be a compact K\"{a}hler manifold with $Ric_{k}(X,\overline{X})\leq-(k+1)\sigma|X|^{2}$, $\sigma\in \mathbb{R}$ and $1<k\leq n$. Then we have
\begin{align}\label{2.3}
(nk+n-k-2)S\cdot h+nRic\leq -n(n+1)(n-1)(k+1)\sigma h \,\,\,.
 \end{align}
\end{lemma}

\noindent {\bf Proof of Proposition \ref{1.2}:} Since the case of $k=1$ is proved by Zhang \cite{Zhang1}, we just consider the case of $1<k\leq n$. Let $\omega\in\alpha$ is a K\"{a}hler metric. For any fixed $k$, we assume $\lambda_{\omega}=\frac{1}{k+1}\mu_{\omega}(k)$, then the $k$-Ricci curvature of $\omega$ satisfies $Ric_{k}^{\omega}\leq (k+1)\lambda_{\omega}$. By the inequality (\ref{2.3}), we have
\begin{align}
(nk+n-k-2)S^{\omega}\cdot \omega +nRic^{\omega}\leq n(n+1)(n-1)(k+1)\lambda_{\omega} \omega \,\,\,,
 \end{align}
and intergrating the above inequality can be obtained
\begin{align}\label{2.4}
(nk+n-k-2)\int_{M}S^{\omega}\omega^{n}&+n\int_{M}Ric^{\omega}\wedge \omega^{n-1} \nonumber \\
&\leq n(n+1)(n-1)(k+1)\lambda_{\omega}\int_{M}\omega^{n} \,\,\,.
 \end{align}
 Note that
 \begin{align}\label{2.5}
 \int_{M}S^{\omega}\omega^{n}=\int_{M}nRic^{\omega}\wedge\omega^{n-1}=-2\pi nc_{1}(K_{M})\cdot \alpha^{n-1}\,\,\,.
 \end{align}
Combining (\ref{2.4}) and (\ref{2.5}), we can get
\begin{align}\label{2.6}
 \lambda_{\omega}\geq \frac{-2\pi c_{1}(K_{M})\cdot \alpha^{n-1}}{(n+1)\alpha^{n}}\,\,\,.
 \end{align}
 Therefore, we obtain
 \begin{align}\label{2.7}
 \mu_{\omega}(k) \geq \frac{-2\pi(k+1) c_{1}(K_{M})\cdot \alpha^{n-1}}{(n+1)\alpha^{n}}\,\,\,,
 \end{align}
 and so
 \begin{align}\label{2.7}
 \mu_{\alpha}(k) \geq \frac{-2\pi(k+1) c_{1}(K_{M})\cdot \alpha^{n-1}}{(n+1)\alpha^{n}}\,\,\,,
 \end{align}
 proving the Proposition \ref{1.2}.            \qed

\section{Proof of Proposition \ref{1.4} and Theorem \ref{1.3} }\label{tk,3}

Let $(M^{n}, \widehat{\omega})$ be a compact K\"{a}hler manifold. The twisted K\"{a}hler-Ricci flow running from $\widehat{\omega}$ satisfies the following equation:
\begin{equation}\label{3.1}
\left\{\begin{split}
&\partial_{t}\omega(t)=-Ric(\omega(t))-\eta \\
&\omega(0)=\widehat{\omega}\,\,\,,
\end{split}\right.
 \end{equation}
 where the $\eta$ is a closed real $(1,1)$ form. It is equivalent to the following Monge-Amp\`{e}re type flow:
 \begin{equation}\label{3.2}
\left\{\begin{split}
&\partial_{t}\varphi=\log \frac{(\widehat{\omega}-tRic(\widehat{\omega})-t\eta+\sqrt{-1}
\partial\overline{\partial}\varphi)^{n}}{\widehat{\omega}^{n}}\,\,\,; \\
&\varphi(0)=0\,\,\,,
\end{split}\right.
 \end{equation}
 Hence, if $\varphi$ is a smooth solution of equation (\ref{3.2}) on $M\times[0,T)$, such that
 \begin{align} \label{3}
 \widehat{\omega}-tRic(\widehat{\omega})-t\eta+\sqrt{-1}
\partial\overline{\partial}\varphi>0 \,\,\,,
 \end{align}
 then $\omega(t)=\widehat{\omega}-tRic(\widehat{\omega})-t\eta+\sqrt{-1}
\partial\overline{\partial}\varphi$ is also a solution of equation (\ref{3.1}). If $\omega(t)$ satisfies
equation (\ref{3.1}), we can define
 \begin{align}
\varphi(t)=\int_{0}^{t}\log\frac{\omega(s)^{n}}{\widehat{\omega}^{n}}ds\,\,\,.
 \end{align}
 We can easily deduce that $\varphi(t)$ satisfies equation (\ref{3.2}).
The solution of the twisted K\"{a}hler-Ricci flow  have a short-time existence (see \cite{GZ}). For the convenience of proof, we need to state three useful lemmas \cite{CLT} which are essentially the same as the case of K\"{a}hler-Ricci flow.
\begin{lemma}\cite{CLT} \label{3.3} Let $\omega(t)$ be a smooth solution to (\ref{3.1}) on $M\times [0, T_{0})$. If there is a positive constant $C>0$ such that
 \begin{align}
C^{-1}\widehat{\omega}\leq \omega(t)\leq C\widehat{\omega} \nonumber
 \end{align}
 on $M\times [0,T_{0})$. Then there is $\varepsilon>0$ such that $\omega(t)$ can be extended to $[0,T_{0}+\varepsilon)$ which satisfies (\ref{3.1}).
\end{lemma}

\begin{lemma} \cite{CLT} \label{3.4}
\begin{equation}
\left\{\begin{split}
& (\frac{\partial}{\partial t}-\Delta_{\omega(t)})\dot{\varphi}=-tr_{\omega}(Ric(\widehat{\omega})+\eta) \,\,\,; \nonumber\\
&(\frac{\partial}{\partial t}-\Delta_{\omega(t)})(t\dot{\varphi}-\varphi-nt)=-tr_{\omega}\widehat{\omega}\,\,\,,\nonumber
\end{split}\right.
 \end{equation}
 where $\omega=\omega(t)$ and $\dot{\varphi}=\partial_{t}\varphi$.
\end{lemma}

\begin{lemma}\cite{CLT} \label{3.5} Let $\omega(t)$ be a smooth solution to (\ref{3.1}) on $M\times [0, T)$.
Then the scalar curvature $S(\omega(t))$ satisfies
\begin{align}
S(\omega(t))+tr_{\omega}\eta\geq-\frac{n}{t+\sigma}  \nonumber
 \end{align}
 on $M\times [0, T)$ where $\sigma>0$ and $inf_{M}(S(\widehat{\omega})+tr_{\widehat{\omega}}\eta)\geq -n\sigma^{-1}$. Moreover,
 \begin{align}
\sup_{M}\log\frac{\det \omega(t)}{\det \widehat{\omega}}=\sup_{M} \dot{\varphi}(\cdot,t)\leq n\log(\frac{t+\sigma}{\sigma}) \,\,\,. \nonumber
 \end{align}
\end{lemma}

We are now ready to prove Proposition \ref{1.4} and Theorem \ref{1.3}.

\noindent {\bf Proof of Proposition \ref{1.4}:}  For any fixed $k$, $1<k<n$, we can assume that $A=\mu_{\widehat{\omega}}(k)>0$. According to the hypothesis of Proposition \ref{1.4}, there exists a positive constant $\delta>0$, such that $\delta[\omega_{0}]+2\pi c_{1}(K_{M})>0$. Hence we can find a smooth function $\upsilon\in C^{\infty}(M)$ which depends on $\widehat{\omega}$, $\widehat{\omega}\in [\omega_{0}]$, such that
 \begin{align}\label{3.6}
\delta\widehat{\omega}-Ric(\widehat{\omega})+\sqrt{-1}\partial\overline{\partial}\upsilon>0 \,\,\,.
 \end{align}

We assume that $\omega(t)$ is the twisted K\"{a}hler-Ricci flow running from $\widehat{\omega}$ with $\eta=\frac{k-1}{2(n-k)}\sqrt{-1}\partial\overline{\partial}\upsilon$. Let $G=tr_{\omega(t)}\widehat{\omega}$.
To simplify notation we write the components of $\omega(t)$ as $g_{i\overline{j}}$ and the components of $\widehat{\omega}$ as $h_{i\overline{j}}$.
Then by the calculation of the parabolic Schwarz Lemma which is a parabolic version of the Schwarz lemma by Yau \cite{Yau2},
we have
 \begin{align}\label{3.7}
(\frac{\partial}{\partial t}-\Delta_{g})\log G\leq \frac{1}{G}g^{i\overline{j}}g^{k\overline{l}}R_{i\overline{j}k\overline{l}}(h)
+\frac{k-1}{2(n-k)}\frac{1}{G}g^{i\overline{l}}g^{k\overline{j}}h_{i\overline{j}}\upsilon_{k\overline{l}} \,\,\,.
 \end{align}
 Applying Lemma \ref{2}, we have
  \begin{align}\label{3.8}
\frac{1}{G}g^{i\overline{j}}g^{k\overline{l}}R_{i\overline{j}k\overline{l}} &\leq \frac{(n-1)A}{2(n-k)}(G+\frac{1}{G}|h|_{g}^{2}) \nonumber \\
&-\frac{(k-1)}{2(n-k)}(tr_{g}Ric)-\frac{(k-1)}{2(n-k)}\frac{1}{G}\langle h,Ric\rangle_{g} \nonumber \\
&=\frac{(n-1)A}{2(n-k)}(G+\frac{1}{G}|h|_{g}^{2})-\frac{(k-1)}{(n-k)}(tr_{g}Ric) \nonumber \\
&+\frac{(k-1)}{2(n-k)}\frac{1}{G}(G\cdot tr_{g}Ric-\langle h,Ric\rangle_{g}) \,\,\,.
 \end{align}
 Choosing local coordinates such that $g_{i\overline{j}}=\delta_{ij}$, $h_{i\overline{j}}=h_{i\overline{i}}\delta_{ij}$, then we also have
 \begin{align}\label{3.9}
G\cdot tr_{g}Ric-\langle h,Ric\rangle_{g}&=\sum_{i}Ric_{i\overline{i}}(\sum_{j}h_{j\overline{j}}-h_{i\overline{i}})\nonumber \\
&=\sum_{i}(Ric_{i\overline{i}}(\sum_{j\neq i}h_{j\overline{j}})) \nonumber \\
&\leq \sum_{i}(\delta h_{i\overline{i}}+\upsilon_{i\overline{i}})(\sum_{j}h_{j\overline{j}}-h_{i\overline{i}}) \nonumber \\
&\leq \delta G^{2}-\delta |h|_{g}^{2}+G\cdot\Delta_{g}\upsilon-\langle \sqrt{-1}\partial\overline{\partial}\upsilon,h\rangle_{g} \,\,\,,
 \end{align}
 where we have used inequality (\ref{3.6}) and $g^{i\overline{l}}g^{k\overline{j}}h_{i\overline{j}}\upsilon_{k\overline{l}}=\langle \sqrt{-1}\partial\overline{\partial}\upsilon,h\rangle_{g}$. Hence, we have
  \begin{align}\label{3.10}
(\frac{\partial}{\partial t}-\Delta_{g})\log G\leq & \frac{(n-1)A}{2(n-k)}(G+\frac{1}{G}|h|_{g}^{2})-
\frac{(k-1)}{(n-k)}(tr_{g}Ric) \nonumber \\
&+\frac{(k-1)\delta}{2(n-k)}G-\frac{(k-1)\delta}{2(n-k)}\frac{1}{G}|h|_{g}^{2}+
\frac{(k-1)}{2(n-k)}\Delta_{g}\upsilon \,\,\,.
 \end{align}
 Since $A>0$, $\delta>0$, $n|h|_{g}^{2}\geq G^{2}$ and $|h|_{g}^{2}\leq G^{2}$, then it is easy to see
   \begin{align}\label{3.11}
(\frac{\partial}{\partial t}-\Delta_{g})\log G\leq & \frac{(n-1)A}{(n-k)}G+\frac{(n-1)(k-1)\delta}{2n(n-k)}G
 \nonumber \\
&+\frac{(k-1)}{2(n-k)}\Delta_{g}\upsilon-\frac{(k-1)}{(n-k)}(tr_{g}Ric)\,\,\,.
 \end{align}
 Combining Lemma \ref{3.4}, we can get
 \begin{align}\label{3.12}
(\frac{\partial}{\partial t}-\Delta_{g})\log G &\leq (\frac{\partial}{\partial t}-\Delta_{g})(-a\xi)+\frac{(k-1)}{2(n-k)}\Delta_{g}\upsilon \nonumber \\
&+\frac{(k-1)}{(n-k)}[(\frac{\partial}{\partial t}-\Delta_{g})\dot{\varphi}+\frac{(k-1)}{2(n-k)}\Delta_{g}\upsilon] \nonumber \\
&=(\frac{\partial}{\partial t}-\Delta_{g})[-a\xi+\frac{(k-1)}{(n-k)}\dot{\varphi}-
\frac{(k-1)}{2(n-k)}\upsilon-\frac{(k-1)^{2}}{2(n-k)^{2}}\upsilon] \,\,\,.
 \end{align}
 where $\xi=t\dot{\varphi}-\varphi-nt$ and $a=\frac{2n(n-1)A+(n-1)(k-1)\delta}{2n(n-k)}$.
 By the maximum principle arguments, we have
 \begin{align}\label{3.13}
\log G\leq C_{1}+(\frac{(k-1)}{(n-k)}-at)\dot{\varphi}+a\varphi+ant \,\,\,.
 \end{align}
 where $C_{1}$ depends on $n$, $k$, $\sup_{M}|\upsilon|$ and $\delta$. This, together with Lemma \ref{3.5}, implies
 \begin{align}\label{3.14}
C^{-1}\widehat{\omega}\leq \omega(t)\leq C\widehat{\omega} \,\,\,.
\end{align}
where $C$ is a positive constant and $t<\min\{T_{\max}, \frac{(k-1)}{a(n-k)}\}$. Combining with Lemma \ref{3.3},
we have
\begin{align}\label{3.15}
T_{max}\geq \frac{(k-1)}{a(n-k)} \nonumber \,\,\,.
\end{align}
Moreover, we conclude that
\begin{equation}\label{3.16}
T_{\max}\geq \frac{2n(k-1)}{(n-1)[2nA+(k-1)\delta]}  \,\,\,. \nonumber
\end{equation}
This completes the proof of Proposition \ref{1.4}.  \qed

\noindent {\bf Proof of Theorem \ref{1.3}:} Since $k=1$ is proved by Zhang \cite{Zhang1}, we only consider the case of $1<k\leq n$.

(1) The case of $k=n$. It means  that $M$ must be of almost nonpositive Ricci curvature. We can choose a family of K\"{a}hler classes $\alpha_{i}$ and a fixed K\"{a}hler metric $\omega_{0}$, such that $\mu_{\alpha_{i}}(n)\alpha_{i}<\frac{1}{i}[\omega_{0}]$, $i=1,2,\cdot\cdot\cdot$.
For any $\alpha_{i}$, we can also find a fixed K\"{a}ler metric $\omega_{\varepsilon_{i}}\in \alpha_{i}$, such that
\begin{equation}
Ric_{\omega_{\varepsilon_{i}}}<\mu_{\alpha_{i}}(n)\omega_{\varepsilon_{i}}+\frac{1}{i}\omega_{0} \,\,\,. \nonumber
\end{equation}
and so
\begin{equation}\label{3.17}
2\pi c_{1}(K_{M})+\mu_{\alpha_{i}}(n)\alpha_{i}+\frac{1}{i}[\omega_{0}]>0\,\,\,.
\end{equation}
We may assume that $\mu_{\alpha_{i}}(n)>0$. If some $\mu_{\alpha_{i}}(n)\leq 0$, we can easily get the result from the Remark \ref{1.31}. Since $0<\mu_{\alpha_{i}}(n)\alpha_{i}<\frac{1}{i}[\omega_{0}]$, we have that
\begin{equation}
c_{1}(K_{M})=\lim_{i\rightarrow \infty}(c_{1}(K_{M})+\frac{\mu_{\alpha_{i}}(n)\alpha_{i}}{2\pi}+\frac{1}{2\pi i}[\omega_{0}]) \,\,\,. \nonumber
\end{equation}
is nef.

(2)  The case of $1<k<n$.

First, we need to prove a fact: if we assume that $\alpha$ is a K\"{a}hler class and $\mu_{\alpha}(k)\geq 0$. We define
\begin{align}\label{3.18}
B_{\alpha}=\inf\{b\in \mathbb{R}| 2\pi c_{1}(K_{M})+b\alpha>0 \} \nonumber \,\,\,.
\end{align}
Then we have $B_{\alpha}\leq \frac{4n}{k-1}\mu_{\alpha}(k)$.

Assume that $B_{\alpha}>\frac{4n}{k-1}\mu_{\alpha}(k)\geq 0$. We shall prove this fact via an argument by contradition. We can choose that $\delta=\frac{3}{2}B_{\alpha}>B_{\alpha}>0$. By the deninition of $B_{\alpha}$, we have
\begin{equation}\label{3.19}
2\pi c_{1}(K_{M})+\delta \alpha>0 \,\,\,.
\end{equation}
For any small $\varepsilon>0$, there is a K\"{a}hler metric $\omega_{\varepsilon}\in \alpha$, such that
\begin{equation}\label{3.20}
\mu_{\omega_{\varepsilon}}(k)<\mu_{\alpha}(k)+\varepsilon \nonumber \,\,\,.
\end{equation}
By the inequality (\ref{3.19}), we can find a smooth real function $\upsilon_{\varepsilon}$ on $M$, such that
\begin{equation}\label{}
-Ric(\omega_{\varepsilon})+\sqrt{-1}\partial\overline{\partial}\upsilon_{\varepsilon}
+\delta\omega_{\varepsilon}>0. \nonumber \,\,\,.
\end{equation}
Let $\eta_{\varepsilon}=\frac{k-1}{2(n-k)}\sqrt{-1}\partial\overline{\partial}\upsilon_{\varepsilon}$.
Applying Proposition \ref{1.4}, we get that  the twisted K\"{a}hler-Ricci flow running from $\omega_{\varepsilon}$,
\begin{equation}
\left\{\begin{split}
&\partial_{t}\omega_{\varepsilon}(t)=-Ric(\omega_{\varepsilon}(t))-\eta_{\varepsilon} \\
&\omega_{\varepsilon}(0)=\omega_{\varepsilon}\,\,\,, \nonumber
\end{split}\right.
 \end{equation}
 exists a smooth solution on $M\times [0, \frac{2n(k-1)}{(n-1)[2n(\mu_{\alpha}(k)+\varepsilon)+(k-1)\delta]})$.
This, together with (\ref{3}),  implies
\begin{align}\label{3.21}
\omega_{\varepsilon}(t)=\omega_{\varepsilon}-tRic(\omega_{\varepsilon})
-t\frac{k-1}{2(n-k)}\sqrt{-1}\partial\overline{\partial}\upsilon_{\varepsilon}+\sqrt{-1}
\partial\overline{\partial}\varphi>0\,\,\,.
\end{align}
where $\varphi$ is a solution of (\ref{3.2}).
From the (\ref{3.21}), it is easy to see
\begin{align}\label{3.22}
\frac{1}{t}\alpha+2\pi c_{1}(K_{M})>0 \nonumber \,\,\,.
\end{align}
where $t\in [0, \frac{2n(k-1)}{(n-1)[2n(\mu_{\alpha}(k)+\varepsilon)+(k-1)\delta]})$. By the deninition of $B_{\alpha}$, we have
\begin{equation}\label{3.23}
B_{\alpha}\leq \frac{(n-1)[2n(\mu_{\alpha}(k)+\varepsilon)+(k-1)\delta]}{2n(k-1)} \,\,\,.
\end{equation}
Since $\varepsilon$ is an arbitrary positive constant, we conclude that
\begin{equation}\label{3.24}
B_{\alpha}\leq \frac{(n-1)[2n\mu_{\alpha}(k)+(k-1)\delta]}{2n(k-1)} \,\,\,. \nonumber
\end{equation}
Combining $B_{\alpha}>\frac{4n}{k-1}\mu_{\alpha}(k)\geq 0$ and $\delta=\frac{3}{2}B_{\alpha}>B_{\alpha}>0$, we have
\begin{align}\label{3.25}
B_{\alpha}&< \frac{(n-1)[2n\cdot\frac{(k-1)B_{\alpha}}{4n}+(k-1)\cdot\frac{3B_{\alpha}}{2}]}{2n(k-1)} \nonumber\\
&=(1-\frac{1}{n})B_{\alpha} \,\,\,. \nonumber
\end{align}
This is a contradiction. So this fact is true.

Now we are ready to prove Theorem \ref{1.3} in case of $1<k<n$. We can choose a family of K\"{a}hler classes $\alpha_{i}$ and a fixed K\"{a}hler metric $\omega_{0}$, such that $\mu_{\alpha_{i}}(k)\alpha_{i}<\frac{1}{i}[\omega_{0}]$, $i=1,2,\cdot\cdot\cdot$. If there is some $\mu_{\alpha_{i_{0}}}(k)=0$, then by the above fact, the $K_{M}$ is nef. If there extists some $\mu_{\alpha_{i_{0}}}(k)<0$,
this theorem has been proved. We assume that $\mu_{\alpha_{i}}(k)>0$ for all $i$. Applying the above fact, we have
\begin{equation}\label{3.26}
\frac{4n}{k-1}\mu_{\alpha_{i}}(k)\geq B_{\alpha_{i}} \,\,\,.
\end{equation}
and so
\begin{equation}\label{3.27}
2\pi c_{1}(K_{M})+\frac{8n}{(k-1)}\mu_{\alpha_{i}}(k)\alpha_{i}>0 \,\,\,,
\end{equation}
where it is a K\"{a}hler classes. Because of $0<\mu_{\alpha_{i}}(k)\alpha_{i}<\frac{1}{i}[\omega_{0}]$, we have
\begin{equation}
c_{1}(K_{M})=\lim_{i\rightarrow \infty}(c_{1}(K_{M})+\frac{4n\mu_{\alpha_{i}}(k)\alpha_{i}}{\pi(k-1)}) \,\,\,. \nonumber
\end{equation}
which implies that $K_{M}$ is nef. We have completed the proof Theorem \ref{1.3}. \qed

\noindent\textbf{Acknowledgement.} The author is grateful to Professor Fangyang Zheng for constant encouragement and support. Many thanks to Yashan Zhang and Man-Chun Lee for useful discussions. The author is also indebted to the work \cite{CLT} and
\cite{Zhang1}.

\end{document}